\documentclass[english,11pt,a4paper]{article}
\usepackage{graphicx}

\usepackage{amsmath}
\usepackage{amsthm}
\usepackage{amscd}
\usepackage{amsfonts}
\usepackage{amssymb}

\usepackage{epsfig}
\newtheorem{definition}{Definition}[section]

\newtheorem{lemma}[definition]{Lemma}

\newtheorem{proposition}[definition]{Proposition}

\newtheorem{remark}[definition]{Remark}

\def\bbb{{\rm l\!I}}
\begin{document}

\title{Quadrature With Respect to Binomial Measures}
\author{F. Calabr\`o         \and
        A. Corbo Esposito\footnote{ DAEIMI \& LAN, Universit\`a di Cassino,
            Via G. Di Biasio 43
            03043 Cassino (FR),
            email \{calabro corbo\} @unicas.it}
}

\date{ArXiv 2, 18/03/08}
\maketitle

\begin{abstract}
This work is devoted to the study of integration with respect to binomial measures. We develop interpolatory quadrature rules and study their
properties. Local error estimates for these rules are derived in a general framework.
\\
Keywords: Numerical quadrature, binomial measures
\\
AMS Subject Classification: 28A25 60G18 65D30 65D32
\end{abstract}

\section{Introduction} \label{intro}

In this paper we develop quadrature formulae for the numerical integration with respect to (w.r.t.) binomial measures. The binomial measure
$\mu_\alpha$, where $0< \alpha < 1$ is a parameter, is a probability measure on an interval of the real line, say $[a,b]$, that is characterized by
the following (self-similar) property: let $I$ a dyadic subinterval of $[a,b]$ and bisect $I$ in the left and right parts $I=I_L \cup I_R$; then
\begin{equation}\label{self_sim}
\mu_\alpha(I_R)=\alpha \mu_\alpha(I) \ .
\end{equation}
When $\alpha=1/2$ we trivially obtain the probability measure proportional to the Lebesgue measure on $[a,b]$. Without loss of generality we will fix
$[a,b]=[0,1]$. The family of measures $\{\mu_\alpha \}_\alpha$ has important features that makes interesting the study of their properties. Firs of
all, it is a family of pairwise mutually singular (see definition 6.7 \cite{rudin}) Borel measures: $\mu_{\alpha_1} \perp \mu_{\alpha_2} $ if
$\alpha_1 \neq \alpha_2 $, thus in particular each $\mu_\alpha$ is singular w.r.t. the Lebeasgue measure $\forall \alpha \ne 0.5$. Moreover,
$\mu_\alpha $ is a continuous measure, i.e. $\mu_\alpha (\{x\})=0\ \forall x\in [0,1] $.
\\
Binomial measures naturally appears in problems related to to the probability theory of sequences of independent trials, and in particular to the so
called Bernoulli process\footnote{Each independent trial has two possible outcomes -success or fail- with fixed probabilities, respectively $\alpha$
and $(1-\alpha)$.}. Let $I_k=\left[\dfrac{j}{2^k},\dfrac{j+1}{2^k}\right[ \subset [0,1[$, where $j,k\in \Bbb N$. Then $\mu_{\alpha}(I_k) $ is exactly
the probability that in the first $k$ trials of the Bernoulli process we have a number of successes equal to the number of $1$s in the binary
expansion\footnote{A detailed study of Hausdorff dimension of sets related to the averages of binary digits is performed in \cite{ACE1,ACE2,ACE3}.}
of $j$.
\\
Also, following the notations  in \cite{Mantica07.12}, the measure $\mu_\alpha$ can be defined as the unique measure that satisfies\footnote{Notice
that the functions $(x/2)$ and $(x/2 + 1/2)$ are usually called Bernoulli shifts.} the following balancing equation:
\begin{equation}\label{shift} 
\int f \, d\mu_\alpha = (1-\alpha) \int f(x/2)\, d\mu_\alpha + \alpha \int f(x/2 + 1/2)\, d\mu_\alpha
\end{equation}
$\forall f\in C^1$. This relation is usefull when multifractal properties are studied,
see \cite{Mandel,Riedi2} for an introduction and \cite{Mantica07.1,Mantica07.2,PeWe2,PeWe} for possible developments.\\[10pt]

The construction of quadrature rules for integration w.r.t. binomial measures is considered for two reasons.
\\
The first is, obviously, for the calculation of the integrals because for these, if  $\alpha \ne 1/2$, we cannot describe the solutions using an
analogous of the fundamental theorem of the integral calculus. Notice, moreover, that the relation \eqref{shift} gives also that the functional
$\displaystyle{L[f]\equiv \int f\,d\mu_\alpha}$ is a "refinable linear functional", as recently defined in \cite{Laurie2006}, with the Bernoulli
shifts as stretch-shift operators and mask $[2(1-\alpha),2\alpha ]$. Thus the calculation of these integrals seams to be an active research problem.
\\
The other reason is that when we exhibit a quadrature rule we are giving in an implicit manner a way, via combination of Dirac delta, for a
decomposition of the measure. Ad example, if the moments are preserved\footnote{This feature is, as we will see, easily related to the degree of
exactness of the quadrature rule.}, the inverse problem is known as moment problem and has been considered also for balanced measures, see
\cite{Abenda,Mantica91,Forte95}. Notice also that the use, for a general measure, of a decomposition that involves the calculation of quadrature
rules based on balanced measures has been explored recently, see \cite{ManticaN2007}.

In this paper we analyze how polynomials can be integrated on dyadic intervals and introduce interpolation based integration rules.

The work  is organized as follows. In the second section we report some analytical background and write some useful technical identities. In the
third section we introduce quadrature with respect to the measures $\mu_\alpha$.  Finally, we list some remarks and possible future work.

\section{Preliminary Results}\label{uno_1}
In this section we present a self-contained introduction of the binomial measures. For the sake of completeness we also list some results that we
will often apply in the sequel. First of all we recall the mean value theorem.
\begin{proposition}[Mean value theorem]\label{mean}
Let $\mu$ be a positive finite measure defined in $[c,d]$, and
consider $\phi\in C^0([c,d])$. Then there $\exists \xi \in [c,d]$:
\begin{equation*}
\phi(\xi)=\dfrac{1}{\mu([c,d])} \int_c^d \phi(x)\, d\mu \ .
\end{equation*}
\end{proposition}
For sake of clearness, and without loosing of generality, we will
restrict ourself to the case of the interval of  integration to be
$I\equiv [0,1]$. An important role will be played by dyadic
intervals $X_j^k \equiv \left[ \dfrac{j}{2^k} , \dfrac{j+1}{2^k}
\right[\ , \quad k\in \Bbb N \,,\ 0\le j< 2^k $. We will call $k$ the
order of the interval.
\\
Let $\mu_{\alpha,k}$ be the probability measure with constant
density on dyadic intervals of order $k$ given by:
\begin{equation*}
f_{\alpha,k}(x)=2^k \alpha^{n(j)} (1-\alpha)^{k-n(j)}\, , \qquad \dfrac{j}{2^k} \le x < \dfrac{j+1}{2^k}\,, \quad j=0\dots 2^k-1
\end{equation*}
where:
\begin{equation*}
n(j)= \# \{ 1\text{s of the binary expansion of }j  \}\ .
\end{equation*}

\begin{figure}
\centering
  \includegraphics[width=0.75\textwidth]{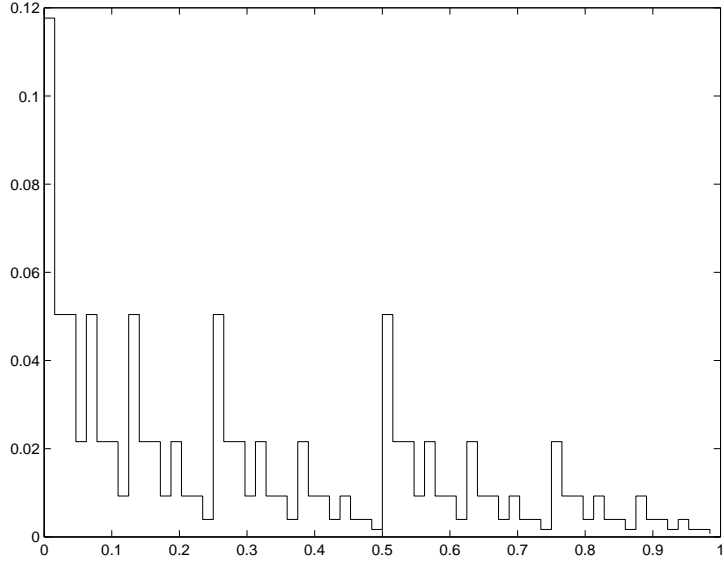}
\caption{The normalization of the density function $2^{-k}f_{\alpha,k}$ with $\alpha=0.3$ and $k=6$}
\label{fig:Funo2}       
\end{figure}

Two useful relations involving this measure are summarized in the
following lemma.
\begin{lemma}\label{lemma1}
Given $h,k \in \Bbb N$, the following hold true:
\begin{gather*}
\mu_{\alpha,k}(X_j^{h})=  \mu_{\alpha,h}(X_j^{h}) \qquad k \ge h\\
\mu_{\alpha,k+h}(X_{j2^h+i}^{k+h})=\mu_{\alpha,k}(X_j^{k})\mu_{\alpha,h}(X_i^{h})\quad \forall i=0\dots 2^h-1, j=0\dots 2^k-1\ .
\end{gather*}
\end{lemma}
\textbf{Proof:} The first relation is immediate, while the second
is easily proved noticing that, due to the fact that since $i$ runs from
$0\dots 2^h-1$ :
\begin{gather*}
n(j2^h+i)=n(j2^h)+n(i)=n(j)+n(i) \Rightarrow\\
\alpha^{n(j2^h+i)} (1-\alpha)^{k+h-n(j2^h+i)} = [\alpha^{n(j)} (1-\alpha)^{k-n(j)}] [ \alpha^{n(i)} (1-\alpha)^{h-n(i)}] \ .
\end{gather*}
\hfill $\Box$
\\
We can now prove the convergence for fixed $\alpha$ of the sequence of measures $\{\mu_{\alpha,k}\}_k$.
\begin{proposition}
The sequence of measures $\{\mu_{\alpha,k}\}_k $ converge in the weak-star sense (see definition 1.58 in \cite{AFP}).
\end{proposition}
\textbf{Proof:} Since the set of probability measures is weak-star
compact, we only need to prove existence of the limit of the
sequence:
\begin{equation*}
\left\{ \int_0^1 \phi(x)\, d\mu_{\alpha,k} \right\}_k \quad \forall \phi \in C^0([0,1])\ .
\end{equation*}
We will prove that this  is a Cauchy sequence. Given $\epsilon
>0$, since $\phi$ is uniformly continuous, there exists   $k_0 \in
\Bbb N$ such that:
\begin{equation*}
\displaystyle{\sup_{|x-y|\le 2^{-k_0} } | \phi(x)-\phi(y)  | < \epsilon}\ .
\end{equation*}
Now, take $k_2 , k_1 \ge k_0$ and consider:
\begin{gather*}
\left| \int_0^1 \phi(x)\, d\mu_{\alpha,k_1}- \int_0^1 \phi(x)\, d\mu_{\alpha,k_2} \right| \le\\
\le \sum_{j=0}^{2^{k_0} -1} \left| \int_{X_j^{k_0} } \phi(x) \, d\mu_{\alpha,k_1}- \int_{X_j^{k_0}} \phi(x)\, d\mu_{\alpha,k_2} \right| =\\
=\sum_{j=0}^{2^{k_0} -1} \left| \phi(\xi_{1,j}) \mu_{\alpha,k_1}(X_j^{k_0})- \phi(\xi_{2,j}) \mu_{\alpha,k_2}(X_j^{k_0}) \right| =\\
= \sum_{j=0}^{2^{k_0} -1} \left| \phi(\xi_{1,j})- \phi(\xi_{2,j})
\right| \mu_{\alpha,k_0}\left(X_j^{k_0}\right) <
\sum_{j=0}^{2^{k_0} -1} \epsilon
\mu_{\alpha,k_0}\left(X_j^{k_0}\right) =  \epsilon \ ,
\end{gather*}
where we have applied the mean value theorem ($\xi_{1,j} ,\xi_{2,j}$ are points in $X_j^{k_0} $) and the first
statement in lemma \ref{lemma1}.\\
The sequence, thus, converges pointwise, as requested. \hfill
$\Box$
\begin{definition}[Binomial measures $\mu_\alpha$]
Fixed $\alpha\in (0,1)$, we will call binomial measure $\mu_\alpha$ the weak-star limit measure of the sequence $\{\mu_{\alpha,k}\}_k$.
\end{definition}
We will denote, as usual, with $L^p_{\mu_\alpha}$ the space of the
$p-$integrable functions with respect to $\mu_\alpha$.
\\
We state now a ``change of variable'' type result for measures $\mu_\alpha$.
\begin{lemma}\label{intqual}
Let $f$ be in $L^1_{\mu_\alpha}$. Then, for each dyadic interval $X_j^k$ we have that:
\begin{equation*}
\int_0^1 f(x)\, d\mu_\alpha = \dfrac{1}{\mu_{\alpha}(X_j^k)} \int_{X_j^k} f(2^k x -j )\, d\mu_\alpha \ .
\end{equation*}
\end{lemma}
\textbf{Proof:} Let us first consider the case where the integrand
$f$ is the characteristic function $\chi_E$ of a measurable set
$E\subset [0,1]$. By lemma \ref{lemma1} the formula is true if $E$
is a dyadic interval, and by summation for a finite union of such
sets. Taking the supremum of such kind of functions we can obtain
the formula for characteristic function of open sets; moreover
taking once again the infimum we can get the formula for any
$\chi_E$ and by linear combination for any simple function.
\\
Let us now take $f\in L^1_{\mu_\alpha}$, $f\ge 0$ (in the general
case we can write $f=f^+ -f^-$ where $f^+$ and $ f^-$ are
respectively the positive and negative parts). Since $f$ is the
pointwise limit of a monotonic sequence of simple functions (see
theorem 1.17 on \cite{rudin}) we obtain the result. \hfill $\Box$
\\[6pt]
Note that lemma  \ref{intqual} cannot be extended to any affine
change of variables, i.e. in general $\int_0^1 f(x)\, d\mu_\alpha
\ne \frac{1}{\mu_\alpha([a,a+ \lambda])} \int_{a}^{a+\lambda}
f(\frac{x-a}{\lambda}) \, d\mu_\alpha$.
\\[6pt]
Now let us see how to calculate the moments. Take $s\in \Bbb N $, we apply lemma \ref{intqual} to write:
\begin{gather}\label{rel1}
\int_0^1 x^s \,d\mu_{\alpha} = \dfrac 1 {1-\alpha } \int_0^{1/2} (2x)^s \,d\mu_{\alpha} = \dfrac {2^s} {1-\alpha } \int_0^{1/2} x^s \,d\mu_{\alpha}
\end{gather}
Now, we can notice that:
\begin{gather*}
\int_0^1 x^s \,d\mu_{\alpha} = \int_0^{1/2} x^s \,d\mu_{\alpha} + \int_{1/2}^1 x^s \,d\mu_{\alpha} =
\end{gather*}
applying another time Lemma \ref{intqual}:
\begin{gather*}
= \int_0^{1/2} x^s \,d\mu_{\alpha} +  \dfrac{\alpha}{1-\alpha} \int_0^{1/2} (x+1/2)^s \,d\mu_{\alpha} = \\
=\int_0^{1/2} x^s \,d\mu_{\alpha} +  \dfrac{\alpha}{1-\alpha} \int_0^{1/2} \left[\sum_{q=0}^s \left( \begin{matrix} s \\ q \end{matrix} \right)   x^{s-q} \dfrac 1 {2^q}  \right] \,d\mu_{\alpha} =\\
= \dfrac 1 {1-\alpha} \int_0^{1/2} x^s \,d\mu_{\alpha} + \dfrac{\alpha}{1-\alpha}  \int_0^{1/2} \left[\sum_{q=1}^s \left( \begin{matrix} s \\ q \end{matrix} \right) x^{s-q} \dfrac 1 {2^q}  \right] \,d\mu_{\alpha}=\\
= \dfrac 1 {1-\alpha} \int_0^{1/2} x^s \,d\mu_{\alpha} + \dfrac{\alpha}{1-\alpha}  \sum_{q=1}^s \left[ \left( \begin{matrix} s \\ q \end{matrix} \right)  \dfrac 1 {2^q} \int_0^{1/2}   x^{s-q} \,d\mu_{\alpha}  \right] =
\end{gather*}
applying the relation \eqref{rel1} in the parenthesis:
\begin{gather*}
= \dfrac 1 {1-\alpha} \int_0^{1/2} x^s \,d\mu_{\alpha} + \dfrac{\alpha}{1-\alpha}  \sum_{q=1}^s \left[ \left( \begin{matrix} s \\ q \end{matrix} \right) \dfrac 1 {2^q} \dfrac{1-\alpha} {2^{s-q}} \int_0^{1} x^{s-q} \,d\mu_{\alpha}  \right]= \\
= \dfrac 1 {1-\alpha} \int_0^{1/2} x^s \,d\mu_{\alpha} + \dfrac{\alpha} {2^s} \sum_{q=1}^s \left[ \left( \begin{matrix} s \\ q \end{matrix} \right)  \int_0^{1}   x^{s-q} \,d\mu_{\alpha}  \right]\ .
\end{gather*}
Summarizing, we have that:
\begin{equation}\label{rel2}
\int_0^1 x^s \,d\mu_{\alpha} = \dfrac 1 {1-\alpha} \int_0^{1/2} x^s \,d\mu_{\alpha} + \dfrac{\alpha} {2^s} \sum_{q=1}^s \left[ \left( \begin{matrix} s \\ q \end{matrix} \right)  \int_0^{1}   x^{s-q} \,d\mu_{\alpha}  \right]
\end{equation}
Now, substituting in the second term of \eqref{rel2} the relation seen in \eqref{rel1} we have that:
\begin{gather}
\int_0^1 x^s \,d\mu_{\alpha} = \dfrac 1 {1-\alpha} \left( \dfrac{1-\alpha}{2^s}  \int_0^1 x^s \,d\mu_{\alpha}  \right) + \dfrac{\alpha} {2^s} \sum_{q=1}^s \left[ \left( \begin{matrix} s \\ q \end{matrix} \right)  \int_0^{1}   x^{s-q} \,d\mu_{\alpha}  \right]  \notag
\end{gather}
and therefore we get the following:
\begin{proposition}\label{momenti_esatti}
Moments of the measures $\mu_\alpha$ are connected by the following recursive relation:
\begin{gather*}
 \int_0^1 x^s \,d\mu_{\alpha} =  \dfrac{\alpha} {2^s-1} \sum_{q=1}^s \left[ \left( \begin{matrix} s \\ q \end{matrix} \right)  \int_0^{1}   x^{s-q} \,d\mu_{\alpha}  \right] \, .
\end{gather*}
\end{proposition}
\hfill $\Box$
\\

We will use piecewise s-polynomial interpolation in the next
chapters. For this reason we explicitly calculate the integrals of polynomials in the
dyadic intervals $X_j^k$ applying lemma \ref{intqual}:
\begin{gather*}
\int_{X_j^k} x^s \,d\mu_\alpha
= \dfrac{\mu(X^k_j)}{(2^k)^s} \sum_{q=0}^s  \left( \begin{matrix} s \\ q \end{matrix} \right) j^q \int_0^1 x^{s-q} \, d\mu_\alpha\ .
\end{gather*}

\section{Quadrature rules}
In this section we want to introduce the numerical integration rules. We will call integration rule a choice of $p+1$
distinct points $ \zeta_q \in [0,1]$ (called nodes) and of values $\beta_q $ (called weights). Let $f(x) \in L^1_{\mu_\alpha}$, we will denote
by\footnote{Note that, beside the integrability condition, we will always  apply quadrature rules to functions with a finite number of
discontinuities and everywhere defined.}:
\begin{equation*}
 \bbb_p (f) \equiv  \sum_{q=0}^p \beta_q  \cdot f(\zeta_q)\, .
\end{equation*}
Notice that, as pointed out in the introduction, this is equivalent to consider as an approximation of the measure $\mu_\alpha$ the following
combination of Dirac delta: $ \mu_\alpha \approx \sum_{q=0}^p \beta_q  \cdot \delta_{\zeta_q} $.
\\
We will call degree of exactness of the formula $\bbb_p$ with respect to $\mu_\alpha$ the greatest positive integer $r$ such that the considered
decomposition maintains the same moments up to order $r$:
\begin{equation}\label{order}
 \int_0^1 x^q \ d\mu_\alpha -\bbb_p(x^q)  =0 \qquad \forall\, q\le r \, ,\ q\in N_0 \ .
\end{equation}
This relation in the case of $r=0$ gives a re-normalization on the weights:
\begin{equation*}
 \sum_{q=0}^p \beta_q  = 1\, .
\end{equation*}
If we consider to replace the function  with a polynomial, we can calculate the integral in an exact manner by means of the seen formulae for
monomials, proposition \ref{momenti_esatti}. If the polynomial is chosen as the one interpolating the function in the nodes, this leads to the so
called interpolation-based integration rules. Fixed the nodes, the weights of the formula  can be calculated integrating, as in the case of the
Lebesgue measure, the so called Lagrange fundamental polynomials (see equation (9.2) in \cite{Qu-Sa-Sa}); in our case this integration is to be made
with respect to the measure $\mu_\alpha$.
\\
A well known theorem, valid for general positive measures (see \cite{GauPol}), states that this rules can give degree of exactness up to $r=2p+1$. It is also well known that, fixed the nodes $\zeta_q, q=0\dots p$ there exists an unique choice of weights that leads to a formula of degree $r\ge p$, and this rule is necessarily interpolation-based.

\begin{table}[t]
\caption{Integration Rules}
\centering
\label{formint2}
    \begin{tabular}{lccc}
\hline\noalign{\smallskip}
 rule & nodes $\zeta_q$ & weights $ \beta_q $ & $\begin{array}{c} \text{degree of}\\ \text{exactness} \end{array}$\\[3pt]
\noalign{\smallskip}
$\Bbb G_0  $ & $\ \alpha $& $\ 1 $ & 1 \\ 
$\Bbb G_1  $ & $\begin{array}{l}   {\small \dfrac {8\alpha+3} {14} - \dfrac  {\sqrt{ -264\alpha^2 +264\alpha +81 }}{42} }   \\[4pt] {\small \dfrac {8\alpha+3} {14} + \dfrac { \sqrt{ -264\alpha^2 +264\alpha +81 }}{42} }  \end{array} $  & $\begin{array}{l}   \dfrac 1 2 - \dfrac{18\alpha-9}{2 \sqrt{ -264\alpha^2 +264\alpha +81 }} \\   \dfrac 1 2 + \dfrac{18\alpha-9}{2 \sqrt{ -264\alpha^2 +264\alpha +81 }} \end{array} $  & 3 \\ 
$\Bbb W_1  $ &  $\begin{array}{c} \alpha - \sqrt{ \dfrac{\alpha(1-\alpha)}{3}  } \\ \alpha + \sqrt{ \dfrac{\alpha(1-\alpha)}{3}  } \end{array}$   &  $\begin{array}{c} 1/2 \\[5pt] 1/2 \end{array}$ & $\begin{array}{cl} 2& if\ \alpha \ne 1/2 \\ 3 & if\ \alpha =1/2 \end{array}$ \\
\noalign{\smallskip}\hline
    \end{tabular}
\end{table}

Optimal rules (from the degree of exactness point of view) are usually referred to as Gauss-Christoffel, see \cite{Gausur} for a complete review in
the general case of a positive measure. We have calculated the first two of them with respect to the measure $\mu_\alpha$ and these are are
summarized in table \ref{formint2}, denoted by $\Bbb G_p$. It can be seen how already $\Bbb G_1 $ is quite laborious to describe \footnote{For the
construction of these rules, the classic procedure is to calculate the nodes as zeros of a polynomial of degree $p+1$ determinate by means of an
orthogonal procedure (weights can be calculated consequently). For a survey about this or other procedures, for example involving the so called
Jacobi matrix, we refer to \cite{Gausur,GauPol}.} for these measures.
\\
If we consider $p+1$ equi-spaced points the corresponding interpolation-based quadrature rules that we obtain are called Newton-Cotes formulae. The
first of these rules are summarized as $\Bbb {NC}_p$ in table \ref{formint3}. We can notice that the weights are negative for some choices of
$\alpha$ at the end-points. Ad example, in the case of $\Bbb {NC}_2$ the weights are positive in the case of $\alpha\in (1/4, 3/4)$, while if $\alpha
< 1/4$ the last weight is negative and in the case of $\alpha>3/4$ the first is negative. Moreover, these formulae are of degree of
exactness\footnote{We have $p+1$ only in the case $\alpha=1/2$ for even $p$.} $p$, that is, obviously, the lowest for interpolation-based rules.
\\
We can moreover fix some parameters and use nonlinear relations \eqref{order} to construct the highest degree interpolatory formula satisfying these
constraints. As example, we can construct a formula with two equal weights and of degree of exactness $2$, see rule $\Bbb W_1$ in table
\ref{formint2}. This formula, although, is meaningful only if $\alpha \in [1/4, 3/4]$ (otherwise nodes are outside of the integration interval) and
coincides with $\Bbb G_1$ in the case $\alpha=1/2$ (for shortness in the sequel we will refer to the case $\alpha=1/2$ as the Lebesgue case).
\\
Another useful example of three point interpolation-based rule  is in table \ref{formint3}, denoted by $\Bbb {GL}_2$. This formula is obtained fixing
the two extreme nodes and considering the only rule of degree of exactness 3 with three nodes. This procedure leads to the so called Gauss Lobatto
formulae. Notice that for this formula the weights are always positive and that in the Lebesgue case $\Bbb {GL}_2 \equiv \Bbb {NC}_2 $.

\begin{table}[t]
\caption{Integration Rules 2}
\centering
\label{formint3}
    \begin{tabular}{lccc} \hline\noalign{\smallskip}
 rule & nodes $\zeta_q$ & weights $ \beta_q $ & $\begin{array}{c} \text{degree of}\\ \text{exactness} \end{array}$ \\[3pt]
\noalign{\smallskip}
$\Bbb {NC}_0  $ & $\dfrac 1  2$ & 1 & $\begin{array}{cl} 0& if\ \alpha \ne 1/2 \\ 1 & if\ \alpha =1/2 \end{array}$ \\ 
$\Bbb {NC}_1  $ &  $\begin{array}{c} 0\\ 1 \end{array}$  & $\left( \begin{array}{c}  1-\alpha\\ \alpha \end{array}\right) $ & 1 \\ 
$\Bbb {NC}_2  $ &$\begin{array}{c} 0 \\ 1/2 \\ 1 \end{array}$& $ \dfrac 1 3 \left( \begin{array}{c}  4 \alpha^2 -  7 \alpha + 3 \\ - 8 \alpha^2 + 8 \alpha \\ 4  \alpha^2 -   \alpha  \end{array} \right) $ & $\begin{array}{cl} 2& if\ \alpha \ne 1/2 \\ 3 & if\ \alpha =1/2 \end{array}$ \\ 
$\Bbb {NC}_3 $ & $\begin{array}{c} 0 \\ 1/ 3 \\  2/ 3 \\ 1 \end{array}$ &  $ \dfrac 1 7 \left( \begin{array}{c}     -9\alpha^3 +24\alpha^2-22\alpha+7 \\  27\alpha^3 -51\alpha^2+24\alpha \\ -27\alpha^3 +30\alpha^2-3\alpha \\ 9\alpha^3 -3\alpha^2+\alpha   \end{array}\right) $ & 3 \\ 
$\Bbb {NC}_4 $ & $\begin{array}{c} 0 \\ 1/ 4 \\  1/2 \\ 3/4 \\ 1 \end{array}$ &  $ \dfrac {1} {315} \left( \begin{array}{c}     256\alpha^4 - 992\alpha^3 + 1572\alpha^2 - 1151\alpha + 315  \\  -32\alpha(32\alpha^3 - 94\alpha^2 +99\alpha-37)  \\ 24\alpha(64\alpha^3 -128\alpha^2+73\alpha-9) \\ -32\alpha(32\alpha^3-34\alpha^2+9\alpha-7)  \\ \alpha(256\alpha^3-32\alpha^2+132\alpha-41) \end{array}\right) $ & $\begin{array}{cl} 4& if\ \alpha \ne 1/2 \\ 5 & if\ \alpha =1/2 \end{array}$ \\ 
$\Bbb {GL}_2  $ &  $\begin{array}{c} 0 \\[5pt] (3\alpha+2)/7 \\[5pt] 1  \end{array}$   &  $ \dfrac 1 3 \left( \begin{array}{c} \dfrac{(\alpha-1)(5\alpha-6)}{(3\alpha+2)} \\[3pt] \dfrac{98\alpha(\alpha - 1)}{(3\alpha + 2)(3\alpha - 5)} \\[3pt] \dfrac{\alpha(5\alpha + 1)}{(5 - 3\alpha)} \end{array} \right) $ & 3 \\
\noalign{\smallskip}\hline
    \end{tabular}
\end{table}

\subsection{A-priori error estimates}\label{par_errestim}
In this section we develop a-priori error estimates for interpolation-based rules. The first that we will see relays on the corresponding formulae
for interpolation errors. Recall that if $\Pi_{f,\zeta_q}(x) $ is the (unique) $p$-polynomial interpolating function $f(x)$ at the nodes $\zeta_q\,,\
q=0,\dots,p$ and $f$ is
sufficiently regular, we have that (see \cite{Qu-Sa-Sa} equation 8.7):  
\begin{equation*}
\begin{array}{l}
f(x)- \Pi_{f,\zeta_q}(x) = \dfrac{f^{(p+1)}(\xi_x )}{(p+1)! }\, \omega_p(x)   \, \text{  for some } \xi_x \in [0,1]\\
\end{array}\,,
\end{equation*}
where $\omega_p$ is the so called nodal polynomial:
\begin{gather*}
\omega_p\equiv \prod_{q=0}^p ( x-\zeta_q )\,.
\end{gather*}
Applying this estimate we obtain:
\begin{gather}
\displaystyle{  \int_0^1 f(x)  \,d\mu_\alpha -\bbb_p (f)  = \int_0^1  \left[f(x)-\Pi_{f,\zeta_q}(x) \right] \,d\mu_\alpha = }\notag \\
\displaystyle{ = \dfrac{1 }{(p+1)!}  \int_0^1    f^{(p+1)}(\xi_x) \omega_p(x)  \, d\mu_\alpha} \label{errcomp1} 
\end{gather}
In the particular case that $f$ is a polynomial of degree $p+k$ with $1\le k \le p+1$ then $f-\Pi_{f,\zeta_q}$ can be factorized as $\omega_p
(x)q(x)$ where $q(x)$ is a polynomial of degree $k-1$. Therefore if $\omega_p$ is $L_{\mu_\alpha}$ orthogonal to all polynomials of degree up to
$k-1$ we obtain that the rule has degree of exactness $p+k$, that is a well known result for quadrature with respect to a general positive measure,
see theorem 2.1 in \cite{GauPol}.
\\
Let us write a simple formula for the error estimate involving
both the derivatives of order $p+1$ and $p+2$. Call $
\omega^+_p(x) = \max\{0, \omega_p(x) \}$, $ \omega^-_p(x) =
\max\{0, - \omega_p(x) \}$, $K_{\alpha}^{+}= \int_0^1
\omega^+_p(x) \,d\mu_\alpha $ and $K_{\alpha}^{-}= \int_0^1
\omega^-_p(x) \,d\mu_\alpha $ then by \eqref{errcomp1} we get:
\begin{gather}
\int_0^1 f(x)  \,d\mu_\alpha -\bbb_p(f)  =\qquad \qquad \qquad \qquad \qquad \qquad \qquad \qquad \qquad \qquad \qquad \qquad \notag \\
\displaystyle{ = \dfrac{1}{(p+1)!}\left[  \int_0^1 f^{(p+1)}(\xi_x) \omega^+_p(x)  \, d\mu_\alpha - \int_0^1 f^{(p+1)}(\xi_x) \omega^-_p(x)  \, d\mu_\alpha \right] } = \tag{\ref{errcomp1}a}\label{errcomp1d} \\
\displaystyle{= \dfrac{1}{(p+1)!}\left[   f^{(p+1)}(\xi_1) \int_0^1 \omega^+_p(x)  \, d\mu_\alpha - f^{(p+1)}(\xi_2) \int_0^1  \omega^-_p(x)  \, d\mu_\alpha \right] } = \tag{\ref{errcomp1}b}\label{errcomp1c}\\
\displaystyle{  = \dfrac{1}{(p+1)!}\left[   K_{\alpha}^{+} [f^{(p+1)}(\xi_1) - f^{(p+1)}(\xi_2) ] - \left[ K_{\alpha}^{-}- K_{\alpha}^{+} \right] f^{(p+1)}(\xi_2)  \right] } = \tag{\ref{errcomp1}c}\label{errcomp1a}\\
\displaystyle{  = \dfrac{1}{(p+1)!}\left[   K_{\alpha}^{+} f^{(p+2)}(\xi_3) (\xi_1-\xi_2)  +  f^{(p+1)}(\xi_2) \int_0^1 \omega_p(x)\,d\mu_\alpha   \right] \qquad \quad} \tag{\ref{errcomp1}d}\label{errcomp1b}
\end{gather}
In \eqref{errcomp1d} we have applied proposition \ref{mean} simply
noticing that $f^{(p+1)}(\xi_x)= \frac{(p+1)![f(x)-\Pi_p (x)]
}{\omega_p(x)}$ can be regarded as a continuous function on
$[0,1]$.
\\
Note that in the same way we could put in evidence $K_{\alpha}^{-} $ instead of $K_{\alpha}^{+} $  equation \eqref{errcomp1a} obtaining an analogous
estimate. This estimate is quite difficult to use, because the $K_{\alpha}$ can be explicitly computed only in few cases (one of these cases can be
seen in the remark \ref{remar2}).
\\[6pt]

\begin{figure}
\begin{center}
\epsfig{file=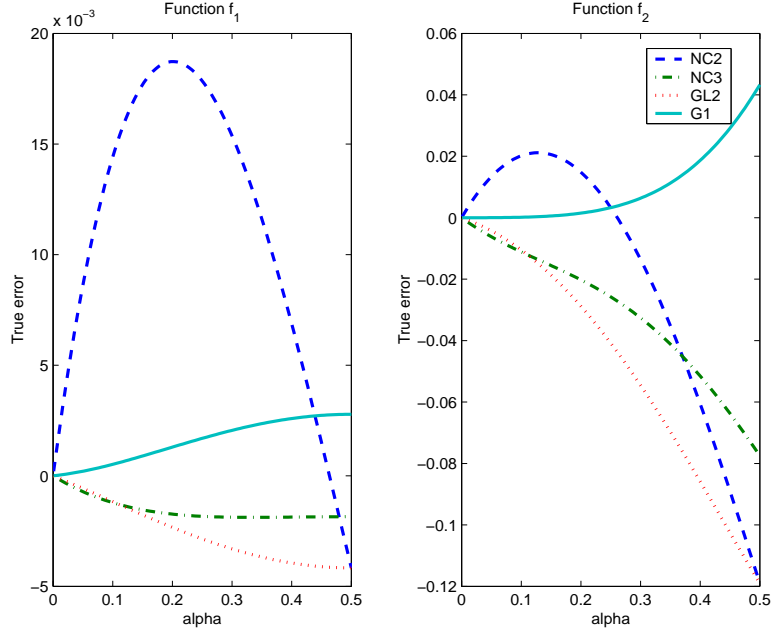,width=0.80\textwidth}
\end{center}\caption{ Plot of the true errors ($\int_0^1 f\,d\mu_\alpha - \bbb(f) $) over the $\alpha$. Test function are in equation \eqref{funz}.}\label{Fig_11}
\end{figure}

Another useful estimate relays on the Taylor expansion of $f$. Consider $\bbb$ be of degree of exactness $r$ and consider $f$ to be $r+1$
times derivable. Take the Taylor expansion of  the function up to the power $r$ in the point $\bar x=0$:
\begin{gather*}
f(x)= P^r_{0}(x) +R^{r+1}_{0}(x) \ \text{ where:} \\
P^r_{0} (x)\equiv \sum_{n=0}^{r} \dfrac{f^{(n)}(0)}{n!} x^n \ ; \quad R^{r+1}_{0} (x)\equiv \dfrac{f^{(r+1)}(\xi_x)}{(r+1)!} x^{r+1} \ .
\end{gather*}
Take now the error of the formula $\bbb$:
\begin{gather}
 \int f(x)\,d\mu_\alpha  - \bbb (f) =\notag \\
 =\int P^r_{0}(x) \,d\mu_\alpha  - \bbb (P^r_{0}) + \int R^{r+1}_{0} (x)\,d\mu_\alpha - \bbb (R^{r+1}_{0}) = \notag \\
= \dfrac{  f^{(r+1)}(\xi) }{(r+1)!}\left[ \int  x^{r+1}  \,d\mu_\alpha - \sum_{q=0} ^N \beta_q \cdot \zeta_q ^{r+1} \right] \label{stim1}
\end{gather}
Notice that we have applied a discrete version of the mean value theorem as in \eqref{errcomp1d}. In this estimate an important role plays the error
committed when calculating the first moment where the formula is not exact. For this reason, in tables \ref{errint1}-\ref{errint2}-\ref{errint3} are
reported errors made by some of the constructed formulae when calculating the moments of the measure. In particular, in each table is fixed $\alpha$.
Errors are quite significant, as it can be seen comparing with the exact values reported in the first line of the tables.
\\
We test some of our formulae also in the case of a polynomial of forth grade and on the $20$-th moment:
\begin{equation}\label{funz}
f_1(x)= \dfrac{5x^4+6x^3-x}{10}\,, \quad f_2=x^{20}
\end{equation}
In figure \ref{Fig_11} we plot the exact errors made from the four formulae $\Bbb {NC}_2$, $\Bbb {NC}_3$, $\Bbb {GL}_2$ and $\Bbb G_1$ when  $\alpha$
varies in $(0,0.5]$.

\begin{remark}[Formulae of degree $p+1$]
Notice that we have constructed, in the the cases of odd nodes, two formulae, namely $ {\Bbb G}_0$ and $ {\Bbb {GL}}_2 $  that generalize the
corresponding equispaced rules and maintain the degree of exactness $p+1$ that is achieved in the Lebesgue case. In the five point case it seems to
us that there is no straight generalization of $\Bbb {NC}_4$ having degree of exactness 5 for all $\alpha$. In particular, if we fix the choice of
the external points to be as in the equispaced case and move the midpoint in function of $\alpha$, the rule exists only for some choices of $\alpha$
(in an interval of $1/2$ approximately of amplitude $0.3185846387$), and the new point is:
\begin{gather*}
\dfrac{320\alpha^3 + 512\alpha^2
-237\alpha-158}{31(64\alpha^2-64\alpha+9)}\ .
\end{gather*}
It is also possible to construct a formulae on 5 points of degree of exactness 5 moving the second and forth point in a rigid linear manner with
respect to $\alpha$, but we obtain a rule that does not coincide with $\Bbb {NC}_4$ in the Lebesgue case. The nodes are:
\begin{gather*}
\left[0, \left(\dfrac{4}{15} - \dfrac{14\sqrt{31} }{465}\right) \alpha - \sqrt{ \dfrac{49\sqrt{31} }{6975}  + \dfrac{1939}{27900} } + \dfrac{7\sqrt{31} }{465} + \dfrac{11}{30} , 1/2,\right. \\
\qquad \left. \left(\dfrac{4}{15} - \dfrac{14\sqrt{31} }{465}\right) \alpha + \sqrt{ \dfrac{49\sqrt{31} }{6975}  + \dfrac{1939}{27900} } + \dfrac{7\sqrt{31} }{465} + \dfrac{11}{30} , 1 \right] \ .
\end{gather*}
We will see in Remark \ref{Extr} how to construct a formula of degree of exactness 5 that coincides with $\Bbb {NC}_4$ in the Lebesgue case but that
uses, in the other cases, six points. \hfill $\Box$
\end{remark}

\begin{table}[h]
\caption{Calculated moments, case $\alpha=0.05$}
\medskip
{\small
\begin{center}
    \leavevmode
    \begin{tabular}{c|c|c|c|c|c} \hline
                         & s= 1                   & 2                  & 3
 & 4                  & 5 \\ \hline
$\int_0^1 x^s\,d\mu_{0.05}$ & $0.0500$  & $0.018\bar 3$ & $0.00860714285714$
 & $0.00448142857143$ & $0.00248679339478$ \\[5pt]
$\Bbb {NC}_0 (x^s)              $ & $0.5000$  & $0.2500$ & $0.12500000000000$
 & $0.06250000000000$ & $0.03125000000000$ \\[5pt]
$\Bbb {NC}_1 (x^s)              $ & $0.0500$  & $0.0500$ & $0.05000000000000$
 & $0.05000000000000$ & $0.05000000000000$ \\[5pt]
$\Bbb {NC}_2 (x^s)              $ & $0.0500$  & $0.018\bar 3$ & $0.00250000000000$
 & $-0.00541666666667$ & $-0.00937500000000$ \\[5pt]
$\Bbb {NC}_3 (x^s)              $ & $0.0500$  & $0.018\bar 3$ & $0.00860714285714$
 & $0.00591798941799$ & $0.00539021164021$ \\[5pt]
$\Bbb {G}_0 (x^s)         $ & $0.0500$  & $0.0025$ & $0.00012500000000$
 & $0.00000625000000$ & $0.00000031250000$ \\[5pt]
$\Bbb {G}_1 (x^s)         $ & $0.0500$  & $0.018\bar 3$ & $0.00860714285714$
 & $0.00407148526077$ & $0.00192634556203$   \\ \hline
    \end{tabular}\label{errint1}
  \end{center}
  }
\end{table}

\begin{table}[h]
\caption{Calculated moments, case $\alpha=0.3$}
\medskip
{\small
\begin{center}
    \leavevmode
    \begin{tabular}{c|c|c|c|c|c} \hline
                         &  s=1                   & 2                  & 3
 & 4                  & 5 \\ \hline
$\int_0^1 x^s\,d\mu_{0.3}$  & $0.3000$  & $0.1600$ & $0.10200000000000$
 & $0.07136000000000$ & $0.05300129032258$ \\[5pt]
$\Bbb {NC}_0 (x^s)              $ & $0.5000$  & $0.2500$ & $0.12500000000000$
 & $0.06250000000000$ & $0.03125000000000$ \\[5pt]
$\Bbb {NC}_1 (x^s)              $ & $0.3000$  & $0.3000$ & $0.30000000000000$
 & $0.30000000000000$ & $0.30000000000000$ \\[5pt]
$\Bbb {NC}_2 (x^s)              $ & $0.3000$  & $0.1600$ & $0.09000000000000$
 & $0.05500000000000$ & $0.03750000000000$ \\[5pt]
$\Bbb {NC}_3 (x^s)              $ & $0.3000$  & $0.1600$ & $0.10200000000000$
 & $0.07511111111111$ & $0.06111111111111$ \\[5pt]
$\Bbb {G}_0 (x^s)         $ & $0.3000$  & $0.0900$ & $0.02700000000000$
 & $0.00810000000000$ & $0.00243000000000$ \\[5pt]
$\Bbb {G}_1 (x^s)         $ & $0.3000$  & $0.1600$ & $0.10200000000000$
 & $0.06725714285714$ & $0.04459836734694$   \\[5pt]
$\Bbb {W}_1 (x^s)        $ & $0.3000$  & $0.1600$ & $0.09000000000000$
 & $0.05080000000000$ & $0.02868000000000$   \\ \hline
    \end{tabular}\label{errint2}
  \end{center}
  }
\end{table}

\begin{table}[h]
\caption{Calculated moments, case $\alpha=0.45$}
\medskip
{\small
\begin{center}
    \leavevmode
    \begin{tabular}{c|c|c|c|c|c} \hline
                        & s=  1                   & 2                  & 3
 & 4                  & 5 \\ \hline
$\int_0^1 x^s\,d\mu_{0.45}$  & $0.4500$  & $0.2850$ & $0.20603571428571$
 & $0.16002428571429$ & $0.13007146313364$ \\[5pt]
$\Bbb {NC}_0 (x^s)              $ & $0.5000$  & $0.2500$ & $0.12500000000000$
 & $0.06250000000000$ & $0.03125000000000$ \\[5pt]
$\Bbb {NC}_1 (x^s)              $ & $0.4500$  & $0.4500$ & $0.45000000000000$
 & $0.45000000000000$ & $0.45000000000000$ \\[5pt]
$\Bbb {NC}_2 (x^s)              $ & $0.4500$  & $0.2850$ & $0.20250000000000$
 & $0.16125000000000$ & $0.14062500000000$ \\[5pt]
$\Bbb {NC}_3 (x^s)              $ & $0.4500$  & $0.2850$ & $0.20603571428571$
 & $0.16373809523810$ & $0.13898809523810$ \\[5pt]
$\Bbb {G}_0 (x^s)         $ & $0.4500$  & $0.2025$ & $0.09112500000000$
 & $0.04100625000000$ & $0.01845281250000$ \\[5pt]
$\Bbb {G}_1 (x^s)         $ & $0.4500$  & $0.2850$ & $0.20603571428571$
 & $0.15456581632653$ & $0.11703565233236$   \\[5pt]
$\Bbb {W}_1 (x^s)         $ & $0.4500$  & $0.2850$ & $0.20250000000000$
 & $0.14805000000000$ & $0.10894500000000$   \\ \hline
    \end{tabular}\label{errint3}
  \end{center}
  }
\end{table}

We conclude this section with two examples on how to apply these estimates \eqref{errcomp1}-\eqref{errcomp1b}, the three point Newton-Cotes rule
$\Bbb {NC}_2$ (where the estimate is quite pessimistic) and the formula $\Bbb {GL}_2$ (where we obtain a good estimate).

\begin{remark}[Error Estimate 1: $\Bbb {NC}_2$]
In the case of the rule with 3 equi-spaced points we can calculate the constants in \eqref{errcomp1b} simply noticing that:
\begin{equation*}
\begin{array}{ll}
\omega_2(x)\ge 0 & if\ x\in [0 , 1/2]\\
\omega_2(x)= -\omega_2(x-1/2) & if\ x\in [1/2 , 1]
\end{array}
\end{equation*}
\begin{equation*}
K^+_\alpha =  \int_0^{1/2} \omega_2(x) \,d\mu_\alpha=  (1-\alpha)\int_0^{1} \omega_2(x/2) \,d\mu_\alpha  = \dfrac{\alpha(1-\alpha)^2 (4-\alpha)}{28}
\end{equation*}
\begin{equation*}
K^-_\alpha = - \int_{1/2}^1 \omega_2(x) \,d\mu_\alpha=  
\dfrac{\alpha}{1-\alpha} \int_{0}^{1/2} \omega_2(x) \,d\mu_\alpha =    \dfrac{\alpha ^2 (1-\alpha)(4-\alpha)}{28}
\end{equation*}
\begin{equation*}
K^+_\alpha -K^-_\alpha = \dfrac{\alpha (1-\alpha)(4-\alpha)(1-2\alpha)  }{28} \ .
\end{equation*}
Taking the absolute value we obtain:
\begin{gather*}
\begin{array}{l}
\displaystyle{\left| \int_0^1 f(x)  \,d\mu_\alpha -\Bbb {NC}_2 f (x) \right| \le}\\[5pt]
\qquad \displaystyle{\le \dfrac{|\alpha (1-\alpha)(4-\alpha)(1-2\alpha) | }{168} \max |f^{\prime \prime \prime}| + }\\[5pt]
\qquad \ \displaystyle{ +\, \dfrac{\alpha(1-\alpha) (4-\alpha)\min(\alpha,1-\alpha)}{168} \max |f^{\prime \prime \prime \prime}| }
\end{array}
\end{gather*}
Notice that this error estimate in the Lebesgue case is pessimistic  because we have $\int_0^1 f(x) \,dx -\Bbb E_2 (f) = -
\dfrac{f^{\prime\prime\prime\prime}(\xi)}{2880}$ for suitable $\xi$. We will obtain the optimal constant in the error estimate in
the next remark when considering the fact that this rule coincides with an Hermite rule.\\
Notice also that this quantity is always null only in the cases $\alpha = 0, 1, 4$ that gives exactly the information that the Newton-Cotes three
points-formulae is an exact formulae only in the dirac cases.
 \hfill $\Box$
\end{remark}

\begin{remark}[Error Estimate 2: $\Bbb {GL}_2$]\label{remar2}
The procedure used above for the estimate of the error is not useful in the case of $\Bbb {GL}_2$ because we know, by construction, that $K_\alpha^+
=K_\alpha^- $, but it is not easy to compute this value, due to the fact that it is given by an integral on a non dyadic subinterval. We can obtain
an error estimate, although, noticing that the integration rule corresponds to integrate the Hermite polynomial with the central node of multiplicity
2. Indeed, see \cite{Qu-Sa-Sa} \P 8.4, the interpolating polynomial in this case can be written as:
\begin{equation*}
\Pi_{f,H}(x) = \sum_{q=0}^2 f(\zeta_q) \lambda_q(x) +   \dfrac{ (x-\zeta_0)( x-\zeta_1 )(x-\zeta_2) }{(\zeta_1-\zeta_0)(\zeta_1-\zeta_2 )  } f^\prime
(\zeta_1)
\end{equation*}
where the functions $\lambda_q(x)$ are the usual Lagrange polynomials. 
Notice that the function that multiplies the $f^\prime(\zeta_1)$
is the nodal polynomial, and this implies that the integral of
this term gives no contribute. By the other side, we can apply the
usual interpolation  error estimate for the Hermite interpolation,
that in our case gives:
\begin{equation*}
\begin{array}{l}
f(x)- \Pi_{f,H}(x) = \dfrac{f^{\prime\prime\prime\prime}(\xi_x ) }{4! }\,  (x-\zeta_0)( x-\zeta_1 )^2 (x-\zeta_2)   \, \text{  for some } \xi_x \in
[0,1]
\end{array} \,,
\end{equation*}
see \cite{Qu-Sa-Sa} \P 8.4. This leads, for the error in the $\Bbb {GL}_2$ case, to:
\begin{gather*}
\int_0^1 f(x)  \,d\mu_\alpha -\Bbb {GL}_2 (f)  = \int_0^1 \dfrac{f^{\prime\prime\prime\prime}(\xi_x ) }{4\,! }\, x(x-1)(x-\zeta_1)^2 \, d\mu_\alpha\
.
\end{gather*}
From this, explicit calculation gives:
\begin{equation*}
-\, \dfrac{2\alpha(17\alpha^3-34\alpha^2+9\alpha+8) }{735}  \dfrac{f^{\prime\prime\prime\prime} (\xi)}{24} \, for\ some \ \xi \in [0,1] \ ,
\end{equation*}
where we have applied proposition \ref{mean} to the function $
f^{\prime\prime\prime\prime}(\xi_x )$ as in \eqref{errcomp1d}.
\hfill $\Box$
\end{remark}

\subsection{Composite rules}
The usual way in which the quadrature rules are used in composite
manner is to introduce a partition of the initial interval and to
consider on each subinterval the integral to be approximated with
a proper quadrature rule. In this paragraph we will see how to do
this in the framework of integration with respect to binomial
measures.
\\
We introduce a partition of the initial interval in $N$ subintervals $J_i=[x_i, x_{i+1}]$ such that $0=x_0 < x_1< \dots< x_N=1  $. We know how to
rescale integrals on dyadic intervals -by means of lemma \ref{intqual}- and for this reason we consider the next definition.
\begin{definition}[Dyadic-regular Partitions]
We will say that the partition $\{J_i\}_{i=0 \dots N-1}$ is dyadic-regular if $\forall i \, \exists j^*$ and $k^*$ s.t. $ J_i=X^{k^*}_{j^*}$. \\
We will say, in particular, that it is dyadic-proper if $N=2^k$ and
$J_i= X^{k}_i \ \forall\, i=0,\dots , 2^k -1$.
\end{definition}

On each subinterval we consider the function to be approximated with a p-polynomial. As seen in the previous chapter, we can use different choices of
interpolating polynomials, leading to different quadrature rules; we will consider that in each subinterval we apply the same interpolatory rule. The
corresponding quadrature formula on the subinterval will be called local quadrature rule. In our notation $ \bbb^N_p (f) $ will indicate that we are
applying the local quadrature rule $ \bbb_p$ on $N$ subintervals to the function $f\in L^1_{\mu_\alpha}$:
\begin{gather*}
\bbb^N_p (f) \equiv \sum_{i=0}^{N-1} \sum_{q=0}^p \beta_q^i f(\zeta_q^i) \,.
\end{gather*}
When we will write $ \bbb^{2^k}_p (f) $ we will consider the dyadic-proper case.
\\
In order to write local quadrature rules we modify nodes and
weights seen in tables \ref{formint2}-\ref{formint3}. In
particular it is easily seen, by lemma \ref{intqual}, that if the
the $p+1$ nodes $\zeta_q^i\,,\ q=1 \dots p+1$  on $J_i \equiv
X_{j^*}^{k^*}$ are taken simply by rescaling them from $[0,1]$ in
the interval $J_i$ $\left(\zeta_q^i \equiv
{(j^*+\zeta_q)}/{(2^{k^*})}\right)$, the corresponding weights are
$ \beta_q^i = \beta_q \mu_\alpha(J_i)$. When we apply the formula
only on one dyadic interval we will use the notation $\bbb_p (f,
J_i)$.
\\[6pt]

\begin{remark}[Extrapolation] \label{Extr}
We have already noticed that both the $\Bbb{NC}_2$ and $\Bbb{GL}_2$ reduce to the 3 points Newton-Cotes formula in the Lebesgue case. In this case
($\alpha=1/2$) we have that the 5 points formula can be written as (see \cite{Esp2003}):
\begin{align}
{\Bbb{NC}}_4  f & =
\dfrac { 16 (\Bbb{NC}_2 (f,X_0^1) + \Bbb{NC}_2(f,X_1^1) ) - \Bbb{NC}_2 (f,X_0^0) }{15} \tag{$\Bbb{NC}^E_2$} \label{forma}\\
&= \dfrac { 16 (\Bbb{GL}_2 (f,X_0^1) + \Bbb{GL}_2(f,X_1^1) ) - \Bbb{GL}_2 (f,X_0^0) }{15} \tag{$\Bbb{GL}^E_2$} \label{formb}
\end{align}
For general $\alpha \ne 1/2$ the rule $\Bbb{NC}^E_2$ uses 5 points and is of degree of exactness 2, while the rule $\Bbb{GL}^E_2$ uses 6 points and
is of degree of exactness 4, thus none of the two are of interpolation type. On the 6 points of the rule $\Bbb{GL}^E_2$ we can construct the rule of
interpolation type that is of order $5$. This formula also coincides with the Newton-Cotes formula $\Bbb{NC}_4$ in the Lebesgue case. We will denote
this formula with ${\Bbb H}_4$. Nodes and weights are:
\begin{gather*}
\begin{array}{rl}
\zeta=&\left[ 0, \dfrac{3\alpha+2}{14}, \dfrac{1}{2}, \dfrac{1}{2}+\dfrac{3\alpha+2}{14}, 1, \dfrac{3\alpha+2}{7} \right]\\
\beta_0=& \dfrac{688\alpha^5-24257\alpha^4 + 59238\alpha^3 -32825\alpha^2-19584\alpha + 16740}{1395(\alpha+3)(3\alpha+2)^2}\\
\beta_1=&\dfrac{224\alpha(1667\alpha^4 -6012\alpha^3 +4855\alpha^2 + 1442\alpha-1952 )}{1395(4-\alpha)(3\alpha-5)(9\alpha^2+12\alpha+4)}\\
\beta_2=&\dfrac{32\alpha(43\alpha^3-86\alpha^2-669\alpha+712)}{1395(3\alpha+2)(5-3\alpha)}\\
\beta_3=& \dfrac{  224\alpha(1667\alpha^4-2323\alpha^3-2523\alpha^2+3395\alpha-216)} {1395(3\alpha-5)(9\alpha^3+18\alpha^2-37\alpha-30)}\\
\beta_4=& \dfrac{\alpha(688\alpha^4+20817\alpha^3-30910\alpha^2-6227\alpha-1108 )}{1395(\alpha-4)(3\alpha-5)^2 }\\
\beta_5=&-\dfrac{98\alpha(2311\alpha^3-4622\alpha^2+1137\alpha+1174)}{1395(81\alpha^4-162\alpha^3-99\alpha^2+180\alpha+100)}
\end{array}
\end{gather*}
Note that the first five nodes are in increasing order, while the last $\zeta_5$ coincides with the midpoint in the Lebesgue case and is between $\zeta_1$ and $\zeta_2$ if $\alpha<1/2$ and between $\zeta_2$ and $\zeta_3$ in the other case.
\hfill $\Box$
\end{remark}

Note that the estimate \eqref{errcomp1b} can be written on a dyadic subinterval by the following:
\begin{gather}
 \int_{X_j^k} f(x)  \,d\mu_\alpha -\bbb_p(f, X_j^k) = \qquad \qquad \qquad \qquad \qquad \qquad \qquad \notag \\
= \dfrac{\mu_\alpha({X_j^k})}{ 2^{k(p+1)} (p+1)!}  \left[  K_{\alpha}^{+}  f^{(p+2)}(\xi^{X_j^k}_3) (\xi^{X_j^k}_1-\xi^{X_j^k}_2) +\right. \qquad \notag \\
\qquad \qquad \qquad \qquad \qquad \left.+ f^{(p+1)}(\xi^{X_j^k}_2) \int_0^1 \omega_p(x)\,d\mu_\alpha \right] ; \label{stima11a}
\end{gather}
By \eqref{stima11a}, taking into account that an analogous equation is true for $K^-_\alpha$, we obtain, for suitable $\eta_1, \eta_2 \in {X_j^k}$:
\begin{gather}
\left| \int_{X_j^k} f(x)  \,d\mu_\alpha -\bbb_p(f, X_j^k) \right| \le \qquad \qquad \qquad \qquad \qquad \qquad \qquad \notag \\
\le \dfrac{\mu_\alpha({X_j^k})}{ 2^{k(p+1)} (p+1)!}  \left[ \dfrac{ \min\{ K_{\alpha}^{+}, K_{\alpha}^{-}  \} }{2^k} \left| f^{(p+2)}(\eta_1)\right| +\right. \qquad \notag \\
\qquad \qquad \qquad \qquad \qquad \left.+ \left|f^{(p+1)}(\eta_2)\right|\left| \int_0^1 \omega_p(x)\,d\mu_\alpha\right|   \right] \,.
\tag{\ref{stima11a}b}\label{stima11b}
\end{gather}
Moreover, we can write in the single dyadic interval also the error estimate \eqref{stim1}:
\begin{gather}
\int_{X_j^k} f(x)\,d\mu_\alpha  - \bbb (f(x), X_j^k)= \notag \\
= \dfrac{  f^{(r+1)}(\xi^{X_j^k}) }{(r+1)!}\left[ \int_{X_j^k}  \left(x-\dfrac{j}{2^k}\right)^{r+1}  \,d\mu_\alpha - \sum_{q=0} ^N  \mu_\alpha(X_j^k) \beta_q \cdot \left( \dfrac{\zeta_q - j}{2^k} -\dfrac{j}{2^k} \right)^{r+1} \right] = \notag \\
=\dfrac{ f^{(r+1)}(\xi^{X_j^k})  \mu_\alpha(X_j^k) }{(r+1)!}  \left( \dfrac{1}{2^k} \right)^{r+1} \left[ \int_0^1 x^{r+1}  \,d\mu_\alpha - \sum_{q=0}
^N  \beta_q \, \zeta_q ^{r+1} \right] \ . \label{stim1b}
\end{gather}

The most simple way to consider the composite integration rules is to consider dyadic proper partitions of increasing order. We report a simple
algorithm that allows us to introduce the definition of order of convergence. We consider  the following iterative procedure:
\begin{gather*}
\textbf{Composite Algorithm} \\
\begin{array}{rl}
\text{\emph{Initialization}} : & \text{put } k=0\\
\text{\emph{Cycle Control}} : & \textbf{while }\text{STOP CRITERION}\\
\text{\emph{Local Quadrature Application}} : &\qquad \text{Compute }  \bbb (f,{X_j^k}) \equiv {\cal I}_j^{(k)} \ \forall j=1\dots 2^k \\
\text{\emph{Bisection}} : &\qquad k=k+1\\\
&\textbf{end while }
\end{array}\\
\left[
\begin{array}{rl}
\text{Output: }& {\cal I} = \sum_{j=1}^{2^k} {\cal I}_j^{(k)} \\
\text{Routine utilized: }& \text{Local quadrature rule } \bbb
\end{array}
\right]
\end{gather*}
We have written the scheme thinking at a procedure that calculates an error estimate and iterates until the result is considered satisfactory: the
stopping criterion has to be chosen to complete the scheme. Usually it is considered in order to satisfy an error requirement:
\begin{equation*}
\left|\int_0^1 f\, d\mu_\alpha - \cal I \right| < tol \ .
\end{equation*}
It can also be  chosen to avoid too many numerical computations. 

\begin{figure}
\begin{center}
\epsfig{file=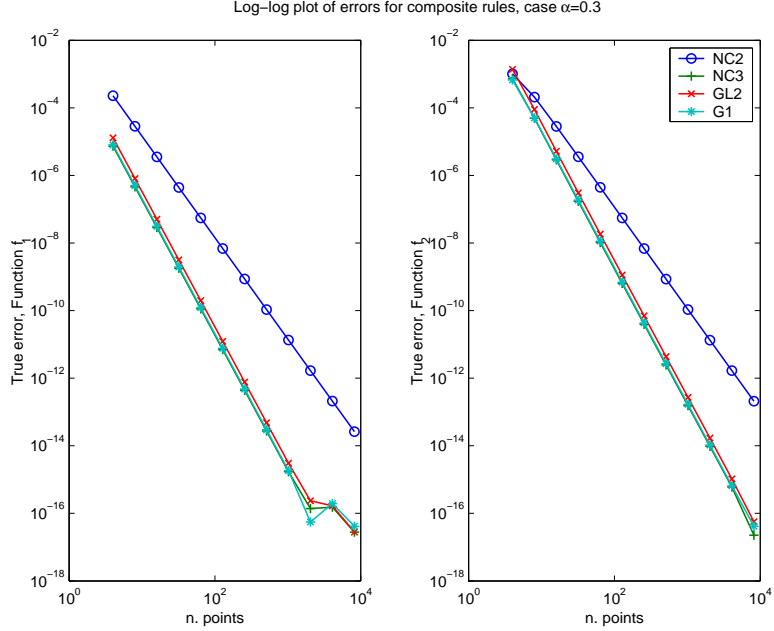,width=0.80\textwidth}
\end{center}\caption{ Plot of the convergence history for some of the integration rules. Test functions are in equation \eqref{funz}.}\label{Fig_12}
\end{figure}

For this algorithm we are interested in convergence properties. We will say that the composite rule converges of order $\gamma$ if
\begin{equation*}
\left| \int_0^1 f(x) d\mu_\alpha - \bbb_p^{2^k} (f)\right| \le K_{f,p} \left(\dfrac{1}{2^k}\right)^{\gamma}\ ,
\end{equation*}
for $f$ sufficiently regular. Applying the Taylor expansion of the function with Peano's remainder \eqref{stim1b} we obtain that the formula has
order of convergence at least equal to the degree of exactness.
\begin{proposition}
Let $\bbb$ be of degree of exactness $r$ and consider $f$ to be $r+1$ times derivable. \\
Then
$\bbb$ has order of convergence at least $r$.
\end{proposition}
\textbf{Proof:} From the error estimate \eqref{stim1b} we have:
\begin{gather*}
\int_{X_j^k} f(x)\,d\mu_\alpha  - \bbb (f, X_j^k)= \notag \\
=\dfrac{ f^{(r+1)}(\xi^{X_j^k})  \mu_\alpha(X_j^k) }{(r+1)!}  \left( \dfrac{1}{2^k} \right)^{r+1} \left[ \int_0^1 x^{r+1}  \,d\mu_\alpha - \sum_{q=0}
^N  \beta_q \, \zeta_q ^{r+1} \right]
\end{gather*}
Now, summing all the subintervals:
\begin{equation}\label{stim2} 
\displaystyle{  \int_0^1 f(x)\,d\mu_\alpha - \bbb_p^{2^k} (f) =  \dfrac{ f^{(r+1)}(\xi) }{(r+1)!} \left( \dfrac{1}{2^k} \right)^{r} k_{\alpha,r, \bbb
}}\, ,  \ for\ some \ \xi \in [0,1]
\end{equation}
where we have called $k_{\alpha,r, \bbb }= \int_0^1 x^{r+1}  \,d\mu_\alpha - \sum_{q=0} ^N  \beta_q \, \zeta_q^{r+1}$. This estimate gives the
requested property.
\hfill $\Box$ \\
In figure \eqref{Fig_12} we plot -for the same two test function seen in the previous section and in the case $\alpha =0.3$- the convergence of the
composite algorithm. From this plot we can see that the order prescribed by the previous theorem is confirmed.

\begin{remark}[Stopping criterion for the non adaptive algorithm]
This convergence property can be used also to write a stopping criterion based on an error estimate, when it is available an estimate of the $p$-th
derivative of the function. The most common strategy, although not very reliable (see \cite{Esp1991,LauNull}), is to estimate the maximum of the
derivative considering the information known from equation \eqref{stim2} at level $k$ and $k+1$. Starting from relation \eqref{stim1b} on a single
interval $X_j^k$ we can write that:
\begin{gather*} \int_{X_j^k} f \, d\mu_\alpha - \bbb (f,X_j^k) = \left(\dfrac 1 {2^k} \right)^{r+1}\; \dfrac{ \mu_\alpha(X_j^k) f^{(r+1)} (\xi^{X_j^k} )}{(r+1)!} k_{\alpha, r,\bbb }
\end{gather*}
and
\begin{gather*}
\begin{array}{l}
\displaystyle{\int_{X_{j}^{k}} f \, d\mu_\alpha - [\bbb (f,X_{2j}^{k+1}) + \bbb (f,X_{2j+1}^{k+1}) ] =} \\
\displaystyle{\qquad = \left(\dfrac 1 {2^{k+1}} \right)^{r+1}\; \dfrac{ \mu_\alpha(X_{2j}^{k+1}) f^{({r+1})} (\xi^{X_{2j}^{k+1}} ) +
\mu_\alpha(X_{2j+1}^{k+1}) f^{({r+1})} (\xi^{X_{2j+1}^{k+1}} ) }{(r+1)!} k_{\alpha,r, \bbb }  \ .}
\end{array}
\end{gather*}
If we consider $f^{({r+1})}$ almost constant in $X_j^k$ to the value $K^{r+1}_{j,k}$ we can consider to approximate this value with:
\begin{gather*}
K^{r+1}_{j,k} \equiv \dfrac{2^{(k+1){(r+1)}} (r+1)! }{(2^{r+1}-1)\mu_\alpha(X_j^k) k_{\alpha,r, \bbb } } \left\{ \bbb (f,X_j^k) - [\bbb
(f,X_{2j}^{k+1}) + \bbb (f,X_{2j+1}^{k+1}) ]\right\}
\end{gather*}
Define, now, $\bar K^{r+1}_k =\max_{j=0,\dots, 2^k-1} K^{r+1}_{j,k} $. With these positions we can consider as stopping criterion in the non-adaptive
composite algorithm:
\begin{equation}
k_{min} \le k \le \min \left\{ \left[ \dfrac 1 {r}  \log_2 \dfrac{k_{\alpha, \bbb } \bar K^{r+1}_{k-1}}{(r+1)! tol} \right]_+ +1 , k_{max} \right\}
\end{equation}
where $[\cdot]_+$ denotes the integer part.

Notice that $k_{min}$ is considered to force the algorithm to do the first computations (ad example in the case of peaked functions) and $k_{max}$ to
avoid too many computations.
\end{remark}
\hfill $\Box$

\section{Final Remarks}
We conclude with some remarks and conclusions.
\\
We develop quadrature formulae for the family of binomial measures. Moreover error estimates, based both on interpolation errors and on Taylor
expansion with Peano's remainder, for such formulae have been established.
\\
We eventually list some possible future developments.
\begin{enumerate}
\item Develop an adaptive algorithm. \item Consider how to calculate numerically rules of higher order. \item To extend the analysis to more general
classes of measures, with special attention to measures given from experimental data.
\end{enumerate}

\bibliographystyle{amsplain}
\bibliography{fract}

\end{document}